\documentclass[11pt]{amsart}
\usepackage{amssymb}

\usepackage[pdftex,colorlinks]{hyperref}
\usepackage[all]{xy}
\usepackage{color}

\theoremstyle{plain} \newtheorem*{theorem}{Theorem}
\theoremstyle{plain} \newtheorem*{corollary}{Corollary}
\theoremstyle{plain} \newtheorem*{lemma}{Lemma}
\theoremstyle{plain} \newtheorem*{proposition}{Proposition}
\theoremstyle{plain} \newtheorem*{definition}{Definition}

\newcommand{\eqq}{==}
\newcommand{\noteqq}{\;\diagup\!\!\!\!\!\!\!{==}}
\newcommand{\existss}{\exists\exists}
\newcommand{\notexistss}{\;\diagup\!\!\!\!\!\!{\existss}}
\newcommand{\inin}{{\in\in}}
\newcommand{\notinin}{\diagup\!\!\!\!\!\!\!\inin}
\newcommand{\impp}{\Longrightarrow}

\newcommand{\id}{\textsf{id}}
\newcommand{\yn}{\textsf{y/n}}
\newcommand{\y}{\textsf{yes}}
\newcommand{\n}{\textsf{no}}
\newcommand{\nt}{\textsf{not}}
\newcommand{\orr}{\textsf{\,or\,}} 
\newcommand{\ad}{\textsf{ and\,}}

\newcommand{\iso}{\textsf{iso}}
\newcommand{\rank}{\textsf{rank}}
\newcommand{\suprank}{\textsf{suprank}}
\newcommand{\minrank}{\textsf{minrank}}

\newcommand{\dom}{\textsf{dom}}
\newcommand{\supt}{\textsf{supt}}

\newcommand{\im}{\textsf{im}}

\newcommand{\card}{\textsf{card}}

\newcommand{\eqv}{\textsf{eqv}}

\newcommand{\morph}{\textsf{morph}}
\newcommand{\obj}{\textsf{obj}}

\newcommand{\PP}{\textsf{P}} 
\newcommand{\bP}{\textsf{bP}}

\newcommand{\sP}{\textsf{sP}}

\newcommand{\WW}{\mathbb{W}}

\newcommand{\lbangle}{\boldsymbol{\langle}}
\newcommand{\rbangle}{\boldsymbol{\rangle}}

\newcommand{\qedms}{\qed\medskip}
\newcommand{\dum}[1]{ \raisebox{1pt}{\framebox[8pt]{$\scriptscriptstyle{#1}$} }}

\begin{document}
\title[A foundation for mathematics]{A foundation for deductive mathematics in the machine age}
\author{Frank Quinn}

\date{November 2026}
\address{Virginia Tech, Blacksburg VA, USA}
\email{fquinn@math.vt.edu}

\begin{abstract} The main outcome is an easy-to-use, modern, fully rigorous set theory. There is a larger environment in which objects have ``elements'' in a weak sense, and the morphisms are the na\"\i ve version of functions. These form a category. The native logic is not 2-valued, and  lacks the 2-valued operations (\(=,\in,\exists, \forall,\nt,\orr)\) of set theory. We then restrict to objects with the minimum needed for sets---an equality pairing and quantification---and see these have all the properties needed for set theory. The weak envioronment enables us to go a step beyond sets: we see that the thing the classical Paradoxes show cannot be a set  lacks only quantification,  and there are several characterizations of it. All this takes less than a dozen pages.

Next we compare this to the leading candidates for foundations in the last century, models of the Zermillo-Fraenkel-Choice axoms. There is a maximal such model in the set theory here, characterized by a property that is assumed in standard mainstream methodology. It follows that only  this maximal model is compatible with our methodology. The ZFC axioms alone do not suffice, as well as being quite a bit more complex. 

In the long run a better foundation will benefit many users of mathematics, especially students. The first applications, however, may be in proof verification. Very recently we were gifted a flood of AI-generated papers on important mathematical topics. Which ones are right? How can we explore them for clues for further progress? The current approach is to ``formalize'' them for input into a proof-checking system (currently LEAN) to see if they can be derived from first principles. The theory here gives a clearer, easier-to-use version of ``first principles'', with a native logic focused more on knowledge than true/false.
 \end{abstract}
 \maketitle
\tableofcontents%

 \section{Introduction} 
 \subsection{Evolution} We locate this work in the historical development. The need to assure reliability led to the ``foundation'' approach: rigorously deduce everything from a  compact core, then lack of contradictions in later development  is evidence for consistency of the core. This, in turn, gives confidence in everything correcty deduced from it. There are still problems with ``correctly deduced'', etc.~but the point is that we didn't \emph{choose} foundations, rigorous deduction, or set theory: this approach became  standard because it  outperformed the alternatives. 
 
Na\"\i ve set theory emerged from the reduction-to-a-core process, but was shown to be inconsistent around 1900. This stalled the development of a foundation, but did not stall the \emph{use} of a foundation. Mainstream mathematicians found that, as long as they stayed away from the known inconsistency (don't say ``set of all sets''), na\"\i ve set theory worked well enough to support more than a century of profound discovery.  Axiomatic set theorists developed ways to work around the problem, but did not resolve it. The development here advances the evolution by  resolving the na\"\i ve inconsistency. The maximality result of Section \ref{sect:wfpairings} also shows we have the entire toolkit:  an axiom not implied by what we already have would make the system inconsistent. 
\subsection{Outline}
The presentation is designed to maximize usefulness, though this puts things slightly out of logical order. We outline the sections, then explain some of the notation. The paper is logically self-contained. The  cores of some of the arguments are quite old, however, so we have included references to \cite{Jech} (\emph{Set theory}, 3rd edition 2002) to connect with this history.
\subsubsection{Primitive set theory} The most useful outcome is the improved set theory. The true primitives and their odd logic are at a lower level, but the oddities only weakly effect the set theory. Accordingly we begin in Section~\ref{sect:sets} with simplified versions ``collections of elements'' and ``functions''. The names should suggest their properties well enough for humans, at least, to work with them. Machines, and people more concerned with category theory, should start with Section \ref{sect:prims+cats}. 
 
`Pre-sets' are collections with a 2-valued `equality' pairing that detects sameness. `Sets' are also required to support quantification. After basic facts about these, well-orders and recursion are described. The universal well-order,  \(\WW\), is  not a set because it does not support quantification. The structure of \(\WW\) and its subobjects is explored in some detail.

\subsubsection{Primitives and categories} Section \ref{sect:prims+cats} begins with a full-precision description of the primitives of the theory. The terminology is changed to better reflect structure: `object generator' replaces `collection of elements' used in \S 1, to reflect similarity to the `object' primitive of category theory, and to inhibit 2-valued logical interpretations. The native logic uses assertions. We might assert ``\(A\) is a set'' because we have either proved this or defined it to be so. However this is not 2-valued: falsity of the assertion does \emph{not} imply ``\(A\) is not a set''. Finally there are four primitive hypotheses. Two are routine, one is a version of the Axiom of Choice, and one concerns Quantification. This last is discussed further in the appendix. 

The usefulness as a setting for category theory is illustrated with a full-precision proof that every ordinary category has a skeleton. 

\subsubsection{Relation to axiomatic set theory} In Section \ref{sect:wfpairings} we show there is a universal well-founded pairing, universal in the sense that any well-founded pairing embeds uniquely as a transitive subobject. In particular, all models of the ZFC  axioms are subobjects of this. Then we see that the universal pairing satisfies the ZFC axioms, so this is the maximal model thereof. The sets here  correspond to sets in this maximal ZFC model.  Section \ref{sect:coherent-lims} identifies properties that characterize the maximal model, and show that  mainstream mathematics implicitly takes place in the set theory defined here.  The additional work (mostly first-order logic) needed for recursion, unions, and some limits to stay within a non-maximal ZFC model is not part of mainstream methodology, so this should not be a complete surprise.
 
 \medskip\par\noindent\textbf{Appendix: omitting ``Quantification''}

The appendix describes the theory developed without the Quantification Hypothesis. The other primitive hypotheses have strong empirical support in traditional work. Traditional logic has  quantification  built in, so does not give clear support for an axiom in a weaker logic. However it turns out that for most practical purposes the theory is unchanged. We therefore assume QH to streamline development and simplify theorem statements.
\subsection{Notation} Some of the differences in this approach require  changes in viewpoint and notation.   We take the opportunty to make other changes to make expressions less ambiguous, better reflect mathematical structure, and easier for both humans and machines to parse. Version of  these are common in computer algebra systems.
\begin{description}
\item[Evaluation] Function evaluation is indicated by square brackets, to distinguish it from parentheses used for grouping: \(f[x]\) is the value of the function \(f\) at \(x\).
\item[Placeholders] In ``define \(f[x]:=x^2\)'' we understand the symbol `\(x\)' as a  placeholder for a value to be specified later.   Here we use `\(\square \)' as a generic placeholder to make this clearer: \(f[\square ]:=\square ^2\) . If the defining expression is short we can use it directly to avoid introducing a symbol for the function. For example   \((\square ^2)[3] =3^2=9\). Finally, in the few cases where we need more than one placeholder we use numbered squares. For instance if we define \(exp[\,\dum1,\dum2]:=\dum1^{\dum2}\), then \(exp[a,b]=a^b\).
\item[Two] Two-valued (Boolean) functions play a major role in set theory. We usually take \(\{\y,\n\}\) (shortened to   `\(\yn\)') as the target of 2-valued functions.  Sometimes we use \(\{1,0\}\) as the target, but we do not shorten it to `\(1/0\)'. We avoid the traditional \(\{\textsf{true, false}\}\) because these terms are used here in non--2-valued ways. 
\end{description}


\section{Primitive set theory}\label{sect:sets}   Set theory is the focus of this section. We go deeper into the roots in the next  section, but a slightly heuristic approach seems workable here because we move quickly to 2-valued logic.  

The primitive objects are called ``\textbf{collections of elements}''. The name should pretty well suggest their nature, but they have no further structure. In particular these do not have an equality pairing that can detect when two elements are the same.  Morphisms of these are \textbf{primitive functions}: \(f\colon A\to B\) is one if it assigns an element of \(B\) to each element of \(A\). These are primitive in the sense that we recognize them by their behavior rather than defining them in terms of something else. Finally, we refer to the whole development as `primitive set theory' to emphasize the nature of the functions.

\subsection{Pre-sets}\label{ssect:pre-sets}

A  \textbf{pre-set} is a collection of elements with a 2-valued equality pairing. Explicitly, \(A\) is a pre-set if there is a function \(=\,\colon A\times A\to \yn\) such that \((a= b)\) implies \(a\) is the same as \(b\), and  the notation `\(=\)' is reserved for this use. ``Same as'' is made more precise in \S\ref{ssect:primitives}.  Equality pairings are unique, so it is sufficient to know when one exists.

\subsubsection{Equivalence classes} 
In brief, the quotient of a collection of elements by an equivalence relation is a  pre-set, and this is the main initial source of pre-sets.

Suppose \(A\) is a collection of elements. An equivalence relation on \(A\) is, as usual, a 2-valued pairing \(\eqv\colon A\times A\to \yn\) that satisfies: \medskip

reflexive: \(\eqv[x,x]=\y\); 

symmetric: \(\eqv[x,y]\implies\eqv[y,x]\); 

and transitive: \((\eqv[x,y]\ad\eqv[y,z])\implies \eqv[x,z]\).\medskip

Given this data we define \textbf{equivalence classes} to be functions \(A\to \yn\) of the form  \(h[\square ]= \eqv[x,\square ]\), for some \(x\in D\).  The collection \(A/\eqv\) of equivalence classes is a pre-set: suppose \(h= \eqv[x,\square ]\) and \(g= \eqv[y,\square ]\), then define \((h= g):= \eqv[x,y]\). The hypothesis that `\(\eqv\)' is an equivalence relation implies that this is a 2-valued pairing, and does not depend on the representatives \(x,y\). Since functions are the same if all values are the same, `\(=\)' implies sameness.

\subsubsection{Example: images}
 Suppose \(f\colon A\to B\) is a function and \(B\) is a pre-set. Then \(\eqv[x,y]:=(f[x]=f[y])\) is an equivalence relation on \(A\), and \(A/\eqv\to B\) is an injective function. This is the \textbf{image} of \(f\),  denoted \(\im[f]\).

\subsection{Sets}
\subsubsection{Quantification}\label{ssect:quant} 

 First, if \(A\) is a collection of elements then the \textbf{powerset} \(\PP[A]\) is the collection of all functions \(A\to \yn\). The \textbf{empty function} is \(\emptyset[\square ]:=\n\), for all \(\square \in A\).
\begin{definition}{\quad}\begin{enumerate}\item A pre-set \(A\) \textbf{supports quantification} if there is a function \(\PP[A]\to\yn\) that detects the empty function. \item A \textbf{set} is a pre-set that supports quantification.
\end{enumerate}
\end{definition} 
Existence of an empty-function detecting function implies the usual quantification expressions (\(\forall a\in A,\text{ or } \exists a\in A\)) are defined and work as usual. 
 
 If the empty function can be detected then all functions can be detected:
\begin{lemma} A pre-set \(A\) supports quantification if and only if \(\PP[A]\) is a pre-set.
\end{lemma} 
Suppose \(\stackrel{\PP}{=}\) is a 2-valued equality pairing on \(\PP[A]\). Then \((\square \stackrel{\PP}{=} \emptyset)\) is a 2-valued function that detects the empty function. 

Conversely, suppose \(\phi\colon \PP[A]\to \yn\) detects \(\emptyset\), ie.~\((\phi[h]= \y)\Leftrightarrow(h== \emptyset)\). Given \(h,g\inin \PP[A]\), the expression \(\square \mapsto (h[\square ]=g[\square ])\) defines a 2-valued function on \(A\), and if this is always `\(\y\)' then the functions are the same. Define `\(\stackrel{\PP}{=}\)' on \(\PP[A]\) by  \[(h\stackrel{\PP}{=} g):= \phi[\square \mapsto\nt[h[\square ] = g[\square ]].\] Then \[(h\stackrel{\PP}{=}g)\Longrightarrow \bigl( \square \in A\Rightarrow h[\square ]= g[\square ]\bigr)\Longrightarrow (h== g).\]\qedms

\subsubsection{ The Quantification Hypothesis}
One of the primitive hypotheses in \ref{ssect:primitives} asserts that if  \(\PP A\) is a pre-set then \(\PP^2A\) is a pre-set. With the definitions above, this is equivalent to: if \(A\) is a set then \(\PP A\) is a set. This is the only hypothesis not already extensively tested in practice. The Appendix outlines the modifications needed to avoid it and explains why, after much experimentation, it seems natural and reasonable to  include it here.

\subsubsection{Subobjects}\label{sssect:subobjects} There are two types of subobjects: images and supports. Images are described above.
If  \(h\colon B\to \yn\) is a 2-valued function then the \textbf{support}, denoted \(\supt[h]\) is the sub-collection on which it has value \(\y\). In the traditional set-builder expression \(\{x\in B\mid h[x]\}\), \(h\) is implicitly assumed to be a 2-valued function, and the expression describes the support of \(h\). 

There is an injection \(\supt[h]\to B\) so supports are images.  The converse is usually true (see below), but fails in some extreme situations. 
 
Standard subset operations may be defined, depending on the type of subobjects. If \(h,g\) are 2-valued functions on \(B\), then intersections and unions are defined by applying logical operations in \(\yn\). For instance \((h\cap g)[\square ]:=(h[\square ]\ad g[\square ])\). A support and an image intersect in an image: if \(f\colon A\to B\) and \(h\) as above are given, then \(\im[f]\cap\supt[h]:=\im[f\upharpoonright  \supt[h\circ f]]\). Unions of images, or of an image and a support, are not defined.  Our first installment in determing when   images are supports is:
\begin{lemma} {\quad}
\begin{enumerate}
\item The image of a set in a pre-set is a support. 
\item A support in a set is a set.
\end{enumerate}
\end{lemma}
For (1), suppose \(f\colon A\to B\) is a function. Then for \(b\in B\), \(f[\square ]=b\) defines a function on \(A\). \(b\) is an element of the image if and only this function is not empty. But on a set we can detect this:  \(b\mapsto ((f[\square ]=b)\stackrel{\PP}{=} \emptyset)\) is a 2-valued function that returns `no' if \(b\) is in the image.

For (2),  if \(h\) is a 2-valued function on \(A\) then a function \(g\) on the support of \(h\) extends to a function on \(A\) by \((h[\square ]\ad g[\square ])\). \(h\) is empty if and only if this extension is empty, and we can detect that.\qedms

Finally, an example of a failure. If \(A\) is a pre-set but not a set, then elements in \(\PP[A]\) are certainly images, but by definition of ``not a set'' these cannot be supports. We expand on this in \ref{sssect:subobjects2}, after we identify a specific pre-set that is not a set. 

\subsubsection{Structure of function domains} \label{ssect:unions}
This includes a strong form of the Union hypothesis of ZFC.
\begin{proposition} \textrm{(Union lemma)} Suppose \(A\) is a collection of elements, \(B\) is a pre-set, and \(f\colon A\to B\) is onto. Then
\begin{enumerate}
 \item \(A\) is a pre-set if and only if all point-inverses \(f^{-1}[b]\) are pre-sets
 
\item\(A\) is a set  if and only if \(B\) and all point-inverses are sets.  
\end{enumerate}
\end{proposition}
Note that since \(B\) is a pre-set, points are supports of 2-valued functions. Composing shows that the preimages in \(A\) are also supports.

Proof of (1): The ``only if'' direction is straightforward, so we focus on ``if''. 
The hypothesis in this direction is that if \(b\in B\)`here is an equality pairing \(\stackrel{b}{=}\) on \(f^{-1}[b]\).  Define \(x\stackrel{A}{=}y\) by: 
\[(x\stackrel{A}{=}y):= \bigl((f[x]\stackrel{B}{=} f[y])\ad (x\stackrel{f[x]}{=} y)\bigr).\]

Proof of (2): First  suppose \(A\) is a set.  Since \(f\) is assumed to be onto, the image of \(f\) is \(B\) and item (1) in the Lemma above shows images of sets are sets. As above point preimages in \(A\) are supports. Item (2) in the Lemma then implies these preimages are sets.

For the other direction, assume \(B\) and the point inverses are sets. Item (1) implies \(A\) is pre-set, so it remains to detect the empty function on it.
Since each \(f^{-1}[b]\) is a set, there are 2-valued functions on \(\PP[f^{-1}[b]]\) that detects emptyness of the restriction. The outputs of these define a 2-valued function on \(B\), and the function on \(A\) is empty if and only if this function is always `yes'. Or equivalently, if the negation is empty. Since \(B\) is a set there is a 2-valued function on \(\PP[B]\) that detects this. Putting these together gives an empty-detecting function on \(\PP[A]\), as required. \qedms\medskip

The name ``Union lemma'' comes from thinking of  the pre-sets \(f^{-1}[b]\) as ``a family of sets, indexed by the set \(B\)''. In these terms the Lemma asserts that the union of such a family is a set. This version is common in the literature.

\subsection{Well-orders and recursion}\label{sect:wellorders} 
We briefly recall the properties of well-orders, describe recursion, and use that to show that there is a universal  well-order.  

\subsubsection{Definitions}\label{ssect:w-oDef}\label{sssect:orders}
Suppose \(A\) is a pre-set.
 A \textbf{linear order}\index{Linear order} on \(A\)  is a 2-valued pairing \((\dum1\leq\dum2)\colon A\times A\to \yn\) that is
\begin{enumerate}\item transitive (\(a\leq b\) and \(b\leq c\) implies \(a\leq c\));
\item any two elements are related: (\((a\leq  b)\orr (b\leq  a)\)); 
 and 
 \item elements related both ways are the same: (\((a\leq  b)\ad (b\leq  a)\Leftrightarrow (a=b)\)).
   \end{enumerate}
  \((A,\leq )\) is a  \textbf{well-order}\index{Well-order} if it is a linear order, and
\begin{enumerate}
 \item if \(a\in A\) then \((\square <a)\) is a set; and
 \item if \(h\colon A\to \yn\) is a  bounded  function with transitive support, then there is \(a\in A\) such that \(h[\square ]=(\square <a)\).
\end{enumerate}
  
\subsubsection*{Notes}
\begin{enumerate}
\item In condition (2) in the definition of well-order:
\begin{enumerate}
 \item  \(h\) `has transitive support' means:  \((h[b]\ad (a\leq b))\) implies \(h[a]\). A function of the form \((\square <a)\) is clearly transitive;
\item  The expression \(h[\square ]=(\square <a)\)  assumes that ``\(\,=\,\)'' is defined for functions on \((\square <a)\). This is why \((\square <a)\) is required to be  a set in item (1).    Note, however,  \(A\) itself is not of the form \((\square <a)\), and  it may not be a set;
\item This  equivalent to: any nonempty support \(B\subset A\) has a \(\leq\)-minimal element. This is the standard definition (cf.~\cite{Jech} definition 2.3), but (2) is slightly easier to work with here.
\end{enumerate}
\item Define a linear \textbf{pre}-order by dropping the requirements that \(A\) be a pre-set and item (3) in the definition of linear order. Then \(
\eqv[a,b]:=(a\leq  b)\ad (b\leq  a)\) is an equivalence relation, and the pre-order induces a linear order on the pre-set of equivalence classes.
\item As usual, linear and well-orders are \emph{hereditary} in the sense that 
if  \((A,\leq )\) is one of these, and  \(h\colon A\to \yn\) is a function, then
 the induced order  on the support of \(h\) is again a linear or well-order.  
\end{enumerate}

\subsubsection{Recursion} \label{ssect:recursion}
Our version is  different from the one given in \cite{Jech}, Theorem 2.15, because we do not find that to be clear. In \S\ref{sect:coherent-lims} we see that being too casual about recursion can conflict with the ZFC axioms. 

First we need a notation for restrictions. Suppose \(f\colon A\to B\) is a partially-defined function and \(D\subset A\) is a set. Then \(f\upharpoonright D\) is the restriction to \(\dom[f]\cap D\). 

\begin{definition}
 Suppose \(A\) is well-ordered and \(B\) is a pre-set. 
\begin{enumerate}
 
\item
 A \textbf{recursion relation} (with values in \(B\)) is a function \(R\colon \dom[R]\to B\) with \(\dom[R]\) some pairs \((a,h)\), where \(a\in A\) and \(h\colon (\square <a) \to B\). 
 \item A partially-defined function \(f\colon A\to B\) is \(R\)-\textbf{recursive} if the domain of \(f\) is transitive, and for every \(a\in\dom[f]\), \(f[a]=R[a,f\upharpoonright (\square <a)]\). 
 \end{enumerate}
\end{definition}
Being \(R\)-recursive is hereditary, in the sense that if \(D\subset \dom[f]\) is transitive then the restriction \(f\upharpoonright D\) is also recursive. 
\begin{proposition}(Recursion) If \((A,\leq), B, R\) are as above, then there is a unique maximal  \(R\)-recursive  function \(f\colon A\to B\). 
\end{proposition}
The proof is essentially the same as the classical one. First, if \(g,h\) are \(R\)-recursive and have the same domain, then they are equal. To see this, suppose not and let \(a\) be the least element on which they differ. Minimality of \(a\) implies the restrictions to \((\square <a)\) are equal. but then 
\[g[a]=R[a,g\upharpoonright (\square <a)]=R[a,h\upharpoonright (\square <a)]=h[a],\]
a contradiction. 

The maximal \(f\) has domain the union \(\cup(\dom[h]\mid h \text{ is }R-\text{recursive})\). If \(x\) is in this union then \(x\in\dom[h]\) for some recursive \(h\). Define \(f[x]:=h[x]\). Uniqueness implies \(f\) is well-defined, and is recursive.\qedms

A recursive \(f\) is maximal if either \(\dom[f]=A\), so not of the form \((\square <a)\), or if \(\dom[f]=(\square <a)\) and \((a,f)\) is not in the domain of \(R\). \medskip

As an aside, a recursion relation can be thought of as a sort of ``vector field'' on  functions with transitive domain, and recursion gives  ``integral curves'' following this vector field. Specifically, a function \(F\) is a `recursive extension' of a function \(f\) defined on \((\square <b)\) if \(F\upharpoonright(\square <b)=f\), and for each \(c>b\), \(F[c]=R[c,F\upharpoonright (\square <c)\, ]\). The recursive extensions of \(f\) trace out the integral curve beginning at \(f\). The proof above shows there is a unique maximal such curve. The traditional result concerns the integral curve beginning with the empty function, with domain the empty set \((\square <\min[A])\).

\subsubsection{Order isomorphisms}\label{sssec:orderisomorphisms}\label{sssect:orderiso}
The first application of recursion is a key result in traditional set theory.
\begin{proposition} Suppose \((A,\leq), (B,\leq)\) are  well-orders. Then there is a unique order-isomorphism from either \(A\) to a  transitive restriction of \(B\), or from a transitive restriction of \(A\) to \(B\). 
\end{proposition}
 This results from the recursion relation \(R[f,a]:=\min[B-f[ (\square <a)]]\). This terminates when either the domain of \(f\) is all of \(A\), so it cannot be of the form \((\square <a)\), or the image is all of \(B\), so \(B-\im[f]\) is empty and there is no minimum. \qedms

\subsection{Universal  well-order}\label{ssect:w-oClassification} 

In this section the well-order \((\WW,\leq)\) is defined, and shown to be universal. This corresponds to the ``ordinal numbers'' of classical set theory, cf.~\cite{Jech}, \S2.

\subsubsection*{Definition} Begin with the collection of all well-ordered \emph{sets}, and define a 2-valued pairing \(\stackrel{w}{\leq}\) on this  by:
\((A,\leq)\stackrel{w}{\leq}(B,\leq)\) if the maximal order-isomorphism described above has domain \(A\). This is a linear pre-order in the sense of Note (2) of \ref{sssect:orders}.  Define  \(\WW\) to be the  pre-set of equivalence classes, with its induced  order. Finally, if  \((A,\leq)\) is a well-ordered set then we denote its equivalence class by \(\lbangle A,\leq\rbangle\) (angle brackets).  

The natural numbers appear here as equivalence classes of \emph{finite} well-ordered sets. Or more invariantly, as the elements less than the first limit point. The smallest element in \(\WW\) is the equivalence class of the empty set, with the vacuous order. It seems appropriate to denote this by ``0'', so the natural numbers here begin with 0 rather than 1.

\subsubsection*{Canonical embeddings}
Recall that if \((A,\leq)\) is a well-order and \(a\in A\),  then \((\square <a)\) with the  induced order is a well-ordered set. It therefore determines an equivalence class \(\lbangle (\square <a),\leq\rbangle\in \WW\). The function  \(\omega[a]:=\lbangle (\square <a),\leq\rbangle\) is the  \textbf{canonical embedding} \(\omega\colon A\to  \WW\).  Since this satisfies the recursion relation  of \ref{sssect:orderiso}, just above, it is the unique order-isomorphism to a transitive restriction of \(\WW\).

\begin{theorem}[Universality of \(\WW\)] {\quad}
\begin{enumerate}\item \((\mathbb{W}, \stackrel{w}{\leq} )\) is a well-order;
\item \(\WW\) does not support quantification, so is not a set;
\item   If \((A,\leq )\) is  a well-order on a set, then the canonical embedding \(\omega\colon A\to\WW\) defined above is an order-isomorphism to the transitive restriction \((\square <\lbangle A,\leq\rbangle)\); if \(A\) is not a set then the image is \(\mathbb{W}\).
\end{enumerate}\end{theorem}
 
\subsubsection*{Proof of the theorem} The first part of  (3) records an observation in the definition of the canonical embedding, above. This implies all the restrictions \(((\square <a),\leq)\) are well-orders, which in turn implies that the order on \(\WW\) is a well-order. Similarly, all the restrictions are sets, so cannot be images of \(A\) if it is not a set.  The alternative (from the recursion relation) is that the image is all of \(\WW\). This gives (1) and completes (3).

We use a form of the Burali-Forti paradox\index{Burali-Forti paradox} to show that  \(\mathbb{W}\) cannot be a set. If it were then \(\WW\) would be a well-ordered set, and would represent an element \(w\in \WW\). The canonical embedding then gives an order-isomorphism of \(\WW\) with \((\square <w)\subset \WW\). The image does not include \(w\), so is not all of \(\WW\).  But the identity is also an order isomorphism, and is onto. This contradicts the uniqueness of order isomorphisms to transitive restrictions.  Therefore \(\WW\) cannot be a set. It is a pre-set, so it must fail to support quantification.
\qedms
 
\subsubsection{Existence of well-orders}\label{ssec:WOexist}
We review a crucial classical consequence (due to Zermillo) of  recursion and the axiom of choice, and get more information about \(\WW\).  
\begin{proposition} Suppose \(A\) is a pre-set. 
\begin{enumerate}\item If \(A\) is  a set then it has a well-order;
\item  \(A\) is  not  a set if and only if there is an injection \(\WW\to A\).
\item There is no 2-valued function on pre-sets that detects which of these alternatives holds.
\end{enumerate}
\end{proposition}
In view of (3), the first two items have to be interpreted as assertions: if we can assert that \(A\) is a set then we can assert that it has a well-order, etc. We include (3) for this perspective, but postpone proof until the next section.\footnote{It also gives a reminder that language alone has no logical force. ``Is \(\square \) a set?'' certainly \emph{sounds} like a 2-valued function, but Mathematics is not impressed.}

Proof: If \(A\) is a set then   \(\sP[A]\) denotes the 2-valued functions on \(A\) whose supports are not all of \(A\). If \(A\) is \emph{not} a set then  \(\sP[A]\) denotes the 2-valued functions whose support is a set. In this case these supports cannot be all of \(A\).
 
Define a \textbf{choice function} for \(A\) to be  \(ch\colon \sP[A]\to A\)  satisfying \(h[ch[h]]=\n\) (ie.~\(ch[h]\) is in the complement of \(h\)). The axiom of Choice implies that any \(A\) has a choice function, as follows:
 Define \(c\colon \sP[A]\times A \to \yn\) by \(c[h,a]:=\nt[h[a]]\). The projection of the support of \(c\) to \(\sP[A]\) is known to be onto, due to the ``not all of \(A\)'' condition. According to Choice, there is a section of this, and sections are exactly  choice functions. 

We now set up for recursion. Fix a choice function and define \(\dom[R]\) to be all pairs \((b,f)\) with \(b\in \WW\) and \(f\colon (\square <b)\to A\). Define \(R\colon \dom[R]\to A\) by : \[R[b,f]:=ch[\im[f\upharpoonright(\square <b)]]).\]
In words, \(f[b]\) is the chosen element in the complement of the image of the restriction \(f\upharpoonright(\square <b)\).
This is clearly a recursion relation. We conclude that there is a maximal partially-defined \(R\)-recursive function \(r\colon \WW\to A\). The definition of \(R\) implies it is an injection, and maximality implies either (1) the domain is  of the form \((\square <a)\) and \(r\) is onto \(A\), or (2) the domain is all of \(\WW\).
\begin{enumerate}
\item If \(A\) is a set then   \(r\) gives a bijection \(r\colon (\square <a)\to A\), for some \(a\in \WW\). The well-order in \(\WW\) restricts to give a well-order on \((\square <a)\), and therefore \(A\) via the bijection.

\item If \(A\) is \emph{not} a set then the domain must be all of \(\WW\), so \(r\colon \WW\to A\) is an injective function.
\end{enumerate}
This completes the proof except for the converse of (2):
\begin{lemma}
 If \(A\) is a set then there is no injective function \(\WW\to A\). 
\end{lemma}
Proof of the lemma: Suppose \(r\colon \WW\to A\) is injective. If \(A\) is a set then \(\PP[A]\) is a set, and by (1) above there is an injection \(\PP[A]\to \WW\). Composing gives an injection \(r\colon\PP[A]\to A\). This contradicts Cantor's theorem, so gives a contradiction. 

For logical completeness we give  details.  \(r\) gives a bijection \(\PP[A]\to \im[r]\), so the inverse \(r^{-1}\colon \im[r]\to \PP[A]\) is, in particular, onto. Next, according to Lemma \ref{sssect:subobjects}, since \(\PP[A]\) is a set, the image of \(r\) is the support of a 2-valued function, say \(h\). Choose an element \(p\in \PP[A]\), and define \(f\colon A\to \PP[A]\) by:
\[f[x]:=\begin{cases} r^{-1}[x] &\text{ if }h[x]=\y\\
p&\text{ if }h[x]=\n,
\end{cases} \]
then this function is also onto.

Now do Cantor. Suppose  \(p\colon A\to \PP[A]\) is a function. Define a 2-valued function \(h\colon A\to \yn\) by \(h[a]:=\nt[(p[a])[a]]\). Then \(h\) is not in the image of \(p\), so \(p\) is not onto. Applying this to  \(f\) gives a contradiction. We conclude there cannot be an injection of \(\WW\) into a set, as claimed.

This completes the proof of the Lemma and therefore the Proposition, except for (3) as noted.\qedms

\subsubsection{Subobjects again} \label{sssect:subobjects2}
The following are corollaries of the Lemma above:
\begin{proposition}{\quad}
\begin{enumerate}
 \item Images in sets are sets; and
 \item Images in \(\WW\) are supports of 2-valued functions.
\end{enumerate}
\end{proposition}
Item (1) now includes images of pre-sets. Suppose   \(B\) is a set and \(A\to B\) is injective. There cannot be an injection \(\WW\to A\), because composing would give an injection into \(B\), contradicting the Lemma. But this means \(A\) must be a set. 

For (2), suppose \(f\colon A\to \WW\) is injective. Then for every \(b\in \WW\) the composition of \((\square <b)\) with \(f\) gives a 2-valued function on \(A\). This is a definition of the intersection, so we denote it by \(f\cap(\square < b)\). \(f\) restricted to this gives an injection into \((\square <b)\). But this last is a set, so by (1) there is a 2-valued function \(h_b\colon \WW\to \yn\) whose support is \(f\cap(\square <b)\). These are coherent in the sense that  if \(b\leq c\) then the restriction of \(h_c\) to \((\square <b)\) is \(h_b\). We can therefore (as in recursion) define \(h\colon \WW\to \yn\) by \(h[b]:=h_c[b]\) for any \(c>b\). This function has the image of \(f\) as its support, as required.
 \qedms
 
The next Lemma uses the above to show some images are not supports, though not in a pre-set. 
\begin{lemma}
 Suppose \(A\) is a pre-set, and there is a nonempty support in \(\PP[A]\) that is a pre-set. Then \(A\) is a set. 
\end{lemma}

Proof: In brief, we can detect elements in a pre-set, so if we can detect a pre-set in \(\PP[A]\) then we can detect some element of \(\PP[A]\). But if we can detect   \(h\in\PP[A]\) then we can detect the empty function, because it is the only one such that both \((\emptyset\orr h)\) and \((\nt[\emptyset]\ad h)\) are  the same as \(h\). Detectability of the empty function is the definition of ``\(A\) supports quantification'', which in turn is the definition of ``\(A\) is a set''. \qedms
 
\begin{corollary} The image of the bounded functions in all 2-valued functions on \(\WW\) is not a support. In other words, the disjoint decomposition of \(\PP[\WW]\) into bounded and unbounded functions is not detected by a 2-valued function.
\end{corollary}
The bounded functions are a pre-set because there is a bijection \(\bP[\WW]\simeq \WW\). This can be seen using the next result, and we find a particularly nice bijection in Section \ref{sect:wfpairings}. Since \(\WW\) is not a set, the Lemma prohibits the existence of a 2-valued function detecting \(\bP[\WW]\subset\PP[\WW]\).\qedms

 \subsection{Identifying \(\WW\)}\label{ssect:idW} This can sometimes be used to determine when a category has isomorphism classes of objects bijective to \(\WW\). The main qualitative point is that in this set theory we can say quite a bit about \(\WW\) beyond just ``not a set''.
\begin{proposition} If \(A\) is a pre-set and \(f\colon A\to \WW\) is a function whose point inverses are sets, then the following are asserted to be equivalent:
\begin{enumerate}
 \item \(A\) is not a set (does not support quantification);
 \item the image of \(f\) is cofinal in \(\WW\); 
 \item there is a bijection \(A\simeq \WW\).
\end{enumerate}
\end{proposition}
Note these  are assertions: if one is known then the others are known, but there is no 2-valued  function that detects  them.

Proof:  Since the point inverses are sets they have well-orders. Choose a well-order  \(\stackrel{a}{\leq}\) on \( f^{-1}[a]\), for every \(a\). Now define an order on \(A\) by: \[(x\leq y):= (f[x]<f[y])\orr (f[x]=f[y] \ad (x\leq_{f[x]}y)).\] This is easily seen to be a well-order, and \(f\) is nondecreasing.  
 The classifying function for this well-order is an order-preserving injection with transitive image, \(\omega\colon A\to \WW\). The Union lemma implies \(A\) is not a set if and only if the image of \(\omega\) is \(\WW\). This gives the implications in the Proposition, except for: if \(\im[f]\subset \WW\)  is cofinal then it, with its induced order, is isomorphic to \(\WW\).
 
 If \(D\subset \WW\) is cofinal then  \(h[a]:=\min[\square \in D \ad \square \geq a]\) defines a surjective function \(\WW\to D\). The point-inverses are bounded, and therefore sets. The domain of \(h\) is not a set so the Union lemma implies \(D\) is not a set. \qedms

\section{Primitives and categories}\label{sect:prims+cats}
Here we describe the primitive foundation, and show how categories fit into this context. Terminology is changed to be more accurately suggestive when fully engaging with the non--2-valued  native logic.  ``collections of elements'' is changed to ``object generators'': ``generator'' because without an `\(\in\)' operation they can't ``collect'' anything; ``object'' to acknowledge the connection to category theory, and because we might not want to call some of these ``elements''. For example we could define a generator whose outputs are  ``feast days in the third year of the reign of Ramses II''. This example also illustrates that we are not trying to define mathematical objects, but rather a \emph{context} in which mathematical objects can be described. Later use of 2-valued functions filters out silly things. 
\subsection{Primitives}\label{ssect:primitives}
 There are three types of primitive ingredients:  objects,  logic, and  hypotheses. 
Some background: Standard practice in mathematics is to define new things in terms of old, and use the definition to infer properties from properties of the old things. But to get started there must be  primitives that are ``recognized'' rather than derived from anything else. As explained in the Introduction, these evolve through trial and error to enable reliable use of mathematical structure. 
The versions here are the latest products of this evoution.

\subsubsection{Primitive objects}\label{ssect:primitiveObj}\index{Primitive!objects}  \medskip
 
\textbf{Object generators} can be thought of as producing things, but without the 2-valued logical structure of traditional sets. There is no equality operator that detects when elements are the same, and there is no membership operator that detects whether or not an element is an output of a given generator. 
 Usage takes the form \(x{\in\in}A\), which we read as ``\(x\) is an output of the generator \(A\)''. 
  
 Syntax for defining generators  takes the form ``\(x\inin A\) means `\dots' ''. For example, the generator whose objects are groups is defined by:\medskip

  \(G\inin\)(groups) means ``\(G\) is a pair \((S,m)\) with \(S\) a set and  \(m\)  a  function \(S\times S\to S\)  that is associative and has a unit and inverses''. \medskip

There is no ``\(\square\in(groups)\)'' function that can identify one of these outputs: we cannot determine if \(G\) is a pair, if the first element is a set, or the second a function. However if we know \(G=(S,m)\) with \(S\) a set and \(m\colon S\times S\to S\)  a function, then there is a 2-valued function of this data that returns `yes' if and only if the function has the desired properties. 

As mentioned above, there are  frivolous generators that are filtered out at a later stage.
 
 \textbf{Morphisms} of  generators are essentially the primitives behind functors of categories, or primitive functions of sets.  ``\(f\colon{\in\in}A\to {\in\in}B\) is a morphism'' means that \(f\) specifies an object \(f[x]{\in\in}B\) for every object \(x{\in\in}A\). 

 Morphisms have some of the structure expected of functors or functions: for instance given morphisms \(\xymatrix{A\ar[r]^f&B}\) and \(\xymatrix{B\ar[r]^g&C}\) then it should be clear that the composition \(\square \mapsto g[f[\square ]]\) qualifies as a morphism.

\subsubsection{Assertion logic}\label{sssect:logic} Primitive logic  provides methods of reasoning with primitive objects and hypotheses. The native logic here is quite a bit weaker than the 2-valued logic used in traditional set theory. We  read traditional 2-valued statements as questions that return `yes' or `no'. For example we read \(\exists (\dots)\) as ``does (\dots) exist?''.  The native logic here uses assertions,  phrased  as ``it is asserted that \dots'', or ``we assert \dots''.  They are denoted by doubling the usual symbols: for example ``\(\existss (\cdots)\) is read as ``we assert \((\cdots)\) exists''. The terminology is a bit awkward, but 2-valued logic is embedded in our language and used in most mathematics, so we need terms that inhibit 2-valued interpretation.
\subsubsection*{Basic assertions} Explanations and clarifications given below.
\begin{enumerate}
\item The basic affirmative assertions are \(\inin,  ==, \existss\);
\item these have negative versions \(\notinin, \noteqq, \notexistss\);
\item  `and'  can be used to combine assertions, but  `or' and `not' require 2-valued logic, so are not available here; and
\item   implication, \(\impp\),  asserts a relationship between assertions.
\end{enumerate} 
We use the common notation \(:=\) for ``is defined by''. Note \((:=)\impp (==)\).
\subsubsection*{Discussion}  {\quad}

\(a\inin A\): This is slightly ambiguous. First, it can mean ``the symbol `\(a\)' is introduced to represent an output of \(A\)''. Second, we can mean ``the previously-defined object \(a\) has been verified to be an output of \(A\)''. In both cases we  assert that it is an output, and usually it doesn't matter how we know this. When it does matter we  rely on the context to determine which is meant. Anyone wanting to dis-ambiguate this could use the notation \( :\!\!\inin\), extending the  ``defined by'' notation, for the first case. \medskip

\(a\eqq b\) is read as  ``\(a\) is asserted to be the same as \(b\)'', and the primitive meaning is that \(a,b\) are symbols representing a \emph{single} output of a generator.  
 This definition makes the usual proof of associativity of composition work for morphisms of generators. Primitive sameness is extremely limited, and going further is mostly restricted to situations where it can be detected with 2-valued logic.\medskip
 
 \emph{Negative versions}: \(a\noteqq b\), for example, is read as ``it is asserted that \(a\) is \emph{not} the same as \(b\)''. These are primitive, and not obtained by somehow negating the affirmative versions.\medskip

\(\impp\): Suppose \(I,J\) are assertions. \(I\impp J\) means  ``if \(I\) can be asserted then \(J\) can also be asserted''. 
In 2-valued logic \(I\impp J\) is a derived relationship defined as \(J\orr \nt[I]\). Here it is primitive because  `\(\orr\)' and `\(\nt\)' are not available. 

 For example, if \(f,g\colon A\to B\) are morphisms of generators then sameness, \(f==g\), means \((a\inin A)\impp (f[a]==g[a])\). In words,   ``if \(a\) can be asserted to be an object in \(A\) then \(f[a]\) can be asserted to be the same as \(g[a]\)''.
Similarly ``\(f\colon A\to B\) can be asserted to be \emph{onto}'' means   \((b\inin B) \impp (\existss a\inin A \mid f[a]==b)\). 
Note that the  `for every' quantifier is subsumed in implication.  \medskip

\emph{Falsehood}:
There are two senses in which assertions can be  false. First, it can be false that we can  assert it.\footnote{
Note that being unable to assert something has no bearing on whether or not it is correct. A philosophical view might be that assertions are more about knowledge than truth or falsehood. } Usually this means the argument that justifies the assertion is missing or erroneous, and this invalidates arguments that depend on the assertion.

The second, and stronger, sense is that we can assert that the thing being asserted is false. This generally means we can assert the negative. For example, if we can assert that the assertion \(a==b\) is false, then we can assert  \(a\noteqq b\).
This is used in contradiction arguments: if all but one input assertion is known to be correct and we can assert the conclusion is false, then the uncertain input can be asserted to be false. We usually interpret this as  asserting  a negative version. However, ``asserting falsehood'' is not an operator on assertions, unlike ``not'' in 2-valued logic. 

 \subsubsection{Primitive hypotheses}\label{ssect:primHyp}
Primitive \textbf{hypotheses} are assertions that we believe are consistent, but cannot justify by reasoning with other primitives. Instead we regard these as experimental hypotheses that have been found to be reliable; see the Discussion.
 \subsubsection*{Hypotheses}
    \begin{description}
\item[Two]
 There is a generator \({\inin}\yn\) such that  \(\y{\inin}\yn\),  \(\n{\inin}\yn\) and if \(a\inin \yn\) then either \(a== \y\) or \(a== \n\), but not both.

\item[Infinity]\index{Axiom!of infinity} The natural numbers support quantification.
\item[Choice] Suppose \(f\colon A\to B\) is a morphism of object generators, and \(f\) is asserted to be onto. Then there is a morphism \(s\colon B\to A\) so that \(f\circ s\) is (asserted to be) the identity. We refer to such \(s\) as \textbf{sections} of \(f\).
\item[Quantification] If a pre-set supports quantification then so does its powerset.
 \end{description}
 
\subsubsection*{Discussion:}\quad

\emph{ About Two}: The force of this hypothesis is that, unlike general  generators, we can tell the objects apart. This axiom is needed for `2-valued functions' to make sense, so is implicit in traditional set theory.\medskip

\emph{About Infinity}: Using primitive objects and the Two and Choice hypotheses, we can construct the natural numbers \(\mathbb{N}\) as a pre-set. However we cannot show that it supports quantification. The Infinity hypothesis asserts this. This  is essentially the same as the ZFC Axiom of Infinity, and is necessary to do much of anything.  \medskip

\emph{About Choice}: The term ``choice'' comes from the idea that if every preimage is nonempty, then we can ``choose'' an element in each preimage to get a morphism \(s\). Note that in general there are no 2-valued functions that determine if a morphism is onto, or if the composition is the identity. These must be assertions, as above. 

The axiom of choice in traditional settings has strong consequences that have been extensively tested for more than a century, and it is now generally accepted. The above  form  extends the well-established version to contexts without quantification.  This extension has implicitly been used in category theory,  again without difficulty. \medskip

\noindent\emph{About  Quantification} Exploration (see the appendix, \S\ref{sect:noQuantification}) shows this hypothesis seems reasonable, does not make much difference, and considerably simplies development and theorem statements.  In the unlikely event that it is shown to be false, it would be straightforward to retroactively complicate statements to accommodate this. 


\subsection{Categories}

 Here we illustrate the convenience of object-generators as a setting for categories. We give a  formally-complete proof  of  the existence of skeleta for ordinary categories. ``Small'' also has a more natural definition here.   

\subsubsection{Definition}
A category  \(\mathcal{C}\) in the terminology here is:
\begin{enumerate}
\item an object generator \(\obj_{\mathcal{C}}\);
\item for every pair \(A,B\inin \obj_{\mathcal{C}}\) a `morphism' object generator \(\morph[A,B]\)
\item for every triple \(A,B,C\inin \obj_{\mathcal{C}}\) a `composition' morphism \[*\colon \morph[A,B]\times\morph[B,C]\to \morph[A,C];\]
\item composition is associative: for every quadruple \(A,B,C,D\inin \obj_{\mathcal{C}}\) the two orders of composition give the same morphism
\[\morph[A,B]\times\morph[B,C]\times \morph[C,D]\to \morph[A,D];\]
\item  for every \(A\inin \mathcal{C}\) there is an `identity' morphism \(\id_A\inin\morph[A,A]\) so that  \(f\inin\morph[A,B]\) implies \(\id_A*f== f== f*\id_B\). 
\end{enumerate}
The `\(*\)' notation for composition reverses the usual order: \(f*g=g\circ f\). We use it here because the expressions are slightly easier to write and parse. 

For example, object generators and their morphisms form a category:  Composition is composition of generator morphisms, and associativity results from our definition of ``same''. 

\subsubsection*{ Domain--range matching error} The composability criterion for morphisms is sometimes stated as: if \(f\in\morph[A,B]\) and  \( g\in \morph[C,D]\), then they are composable if \(B= C\).  This only makes sense if the collection of objects is a pre-set, which is rarely the case in nature. It is true for skeleta (see below), but these are not particularly natural.  The assertion form \(B== C\) is not 2-valued so does not provide a test for composability. 

\subsubsection{Isomorphism classes}\label{ssect:isoclasses}

Suppose  \(X,Y{\inin}\mathcal A\) are objects in a category. An \textbf{isomorphism} is, as usual, a morphism \(i\inin\morph[X,Y]\) such that there exists a morphism \(j\inin\morph[Y,X]\) so that the compositions \(i*j\) and \(j*i\) are identities.  If \(i\) is an isomorphism then, again as usual, the inverse \(j\) is unique, and is also an isomorphism. Composition of isomorphisms give isomorphisms. The isomorphisms in a category thus define a subcategory, which we denote by \(\iso[\mathcal A]\).

A category is \textbf{ordinary} if the morphism generators are sets. In this case there are 2-valued functions that detect the identity morphisms, so ``are  \(i*j\) and \(j*i\) identity morphisms?'' defines a 2-valued function on \(\morph[X,Y]\times\morph[Y,X]\). Therefore  there is a 2-valued function of \(X,Y\) that implements ``is there an isomorphism \(X\to Y\)?''s. This gives an equivalence relation on objects, and the equivalence classes are, as usual, a pre-set. We denote this by  \(\obj \mathcal A/\iso\). There is a quotient morphism (of object generators) \(q\colon \obj\mathcal{A}\to \obj\mathcal{A}/\iso\). 

Since functors preserve identity morphisms, they also preserve isomorphisms. A functor \(F\colon \mathcal{A}\to \mathcal{B}\) therefore induces a function \(F\colon \obj\mathcal{A}/\iso\to \obj\mathcal{B}/\iso\). In other words \(\obj\square /\iso\) is a functor, from ordinary categories and functors, to pre-sets and functions.  

\subsubsection{Skeleta}\label{ssect:skeleta}
An ordinary category \(\mathcal{S}\) is \textbf{skeletal} if  isomorphic objects are equal. 
This is equivalent to: the quotient function, from objects to isomorphism classes, is a bijection. 

A  \textbf{skeleton}\index{Skeleton} of a category  \(\mathcal{A}\) is a skeletal category and a functor \(s\colon\mathcal{S}\to \mathcal{A}\) that  induces a bijection on isomorphism classes of objects, and is a bijection on morphism sets.

\begin{theorem}[Existence of skeleta]{\quad}
\begin{enumerate}\item Every ordinary category has a skeleton; 
\item the inclusion of a skeleton is an equivalence of categories. More precisely, if \(s\colon \mathcal{S}\to \mathcal{A}\) is the inclusion of a skeleton then there is a functor \(q\colon \mathcal{A}\to \mathcal{S}\) so that \(q\circ s\) is the identity of \(\mathcal{S}\) and \(s\circ q\) is naturally equivalent to the identity of \(\mathcal{A}\).
\end{enumerate}
\end{theorem}
\begin{corollary}(Uniqueness) If \(s_n\colon \mathcal{S}_n\to \mathcal{A}\) for \(n=1,2\) are skeleta, then there is an isomorphism of categories \(T\colon  \mathcal{S}_1\to  \mathcal{S}_2\) so that \(s_2\circ T\) is naturally equivalent to \(s_1\).
\end{corollary}
For the Corollary, let \(q_2\colon \mathcal{A} \to \mathcal{S}_2\) be the natural inverse in item (2) of the proposition. Then \(T=p_2\circ s_1\) is the desired isomorphism.\qedms

We give details for people---and machines---who want to see how the primitive logic and hypotheses mesh to enable full-precision proofs.  Others should skip this. 

For the first step note that the quotient \(q\colon \obj\mathcal{A}\to \obj\mathcal{A}/\iso\) is known to be a surjective morphism of generators. Choice therefore asserts that there is a section, \(h\colon  \obj\mathcal{A}/\iso\to \obj\mathcal{A}\). Define a category \(\mathcal{S}\) with objects the equivalence classes \(\obj\mathcal{A}/\iso\) and morphisms \(\morph_{\mathcal{S}}[x,y]:=\morph_{\mathcal{A}}[h[x],h[y]]\). It should be clear that \(\mathcal{S}\) is skeletal. We get a functor \(s\colon \mathcal{S}\to \mathcal{A}\) by the section \(h\) on objects, and the identity on morphisms. It should be clear that this is a skeleton.

The next step is to extend the quotient function \(p\colon \obj\mathcal{A}\to \obj\mathcal{A}/\iso=\obj\mathcal{S}\) to a functor. 
To begin, consider the object generator with objects \((a,\theta)\), where \(a\) is an object in \(\mathcal{A}\) and \(\theta\) is an isomorphism \(a\simeq h[q[a]]\). The forgetful morphism from this to \( \obj\mathcal{A}\) is onto, so again we can chose a section. Denote this by \(a\mapsto (a, \theta_1[a])\). We tidy this up a bit.  Define \(\theta_2\) by 
\[\xymatrix{ a\ar[r]^{\theta[a]}& h[a]\ar[r]^{\theta[h[a]]^{-1}}&h[a].
}\]
Then \(\theta_2\) has the benefit that it takes an image object \(h[b]\) to \((h[b],\id)\). Note that we cannot get \(\theta_2\) by saying ``for objects in the image of \(h\), define \(\theta_2[h[a]]:=(h[a],\id)\) and extend randomly to other objects''. For this to be valid we would need a 2-valued function on \( \obj\mathcal(A)\) that detects the image of \(h\), and these may not exist. 

Now extend \(q\) to morphisms by: for \(f\colon a\to b\), \(q[f]\) is the composition
\[\xymatrix{h[q[a]]\ar[r]^{\quad\theta_2[a]^{-1}}&a\ar[r]^{f}&b\ar[r]^{\theta_2[b]\quad}&h[q[b]]}\]
This is a morphism in \(\mathcal{S}\) because \(h\) is the identity on morphisms. 

The proof is completed by the routine observations that \(q\) is a functor; \(q\circ h=\id[\mathcal{S}]\); and \(\theta_2\) is a natural equivalence \(h\circ q \simeq \id[\mathcal{A}]\). \qedms


\section{Well-founded pairings}\label{sect:wfpairings}
This is the first part of the comparison of the set theory here with classical axiomatic set theory as exposed in \cite{Jech}. The key to the comparison is the ``Well-founded''  axiom,  proposed by von Neumann in 1925 to avoid self-reference and circularity.

A pairing  \(\lambda\colon U\times U\to \yn\) on a pre-set is \textbf{well-founded} if 
\begin{enumerate}
 \item \(\lambda[\square ,a]\) is a set, for all \(a\in U\), and
 \item  the complement  (in \(U\)) of a set has a \(\lambda\)-minimal element. 
\end{enumerate}
 Here  \(a\in S\) is said to be \(\lambda\)-\textbf{minimal} in \(S\) if \(\lambda[x,a]\Longrightarrow x\notin S\).

The \textbf{adjoint} of such a pairing  is the function \(\lambda^{adj}\colon U\to \bP[U]\) defined by \(\lambda^{adj}[a]:=\lambda[\square ,a]\). Technically this is the adjoint of the opposite, but since this is the only use we make of `adjoint', simplifying the notation seems harmless. `\(\bP\)' refers to 2-valued functions whose rank (described below) is bounded; boundedness is needed to ensure the range of the adjoint is a pre-set.

A \textbf{morphism} of pairings \((U,\lambda)\to (V,\mu)\) is a function \(f\colon U\to V\) that satisfies \(\lambda[a,b]=\mu[f[a],f[b]]\). In adjoint terms, this is commutativity of the diagram 
\[\xymatrix{U\ar[r]^f\ar[d]^{\lambda^{adj}}&V\ar[d]^{\mu^{adj}}\\
\bP[U]&\bP[V]\ar[l]_{f^*}}\]
\subsection{Main theorem}
\begin{theorem} {\quad} 
\begin{enumerate}
\item There is a well-founded pairing \(\epsilon\colon \WW\times\WW\to \yn\) whose adjoint gives a bijection \(\WW\to \bP[\WW]\);
\item Pairings as in (1) are universal for well-founded pairings  \((U,\lambda)\) with injective adjoint, in the sense that there is a unique injection \(U\to \WW\) that is a morphism of pairings and has transitive image;
\item Pairings satisfying the ZFC axioms are well-founded in this sense (possibly after truncation so the range of `rank' is \(\WW\); see \ref{sssect:setsRsets} in the next section);
\item Pairings as in (1) satisfy the ZFC-1 axioms, where `-1' means ``omit first-order logic requirements''; and finally
\item The set theory associated to pairings as in (1) is the primitive set theory of \S 2, restricted to subsets of \(\WW\).
\end{enumerate}
\end{theorem}
Consequences: well-founded pairings on \(\WW\) with bijective adjoint  (as in (1)) are isomorphic, and give a maximal model for the ZFC axioms. The universality  (2) is  a version of Mostowski collapsing (\cite{Jech}, \S6.15). The injective-adjoint condition is traditionally called ``extensional'', and corresponds to the set-theory axiom that sets are determined by their elements. 

\subsubsection{Proof of (5)} This one is easy so we do it immediately. Surjectivity of the adjoint means every primitive-theory subset is a set in the associated model of the ZFC axioms. A primitive function \(A\to B\) is a function in the model if and only if its graph in \(A\times B\) is a set in the model. But the graph is a primitive-theory subset, and (again) surjectivity of the adjoint implies this is in the model. \qed
\subsection{Rank functions} 
\begin{lemma} A pairing \(\lambda\) is well-founded if and only if 
  there is a function \(\rank\colon U\to \WW\) so that for every \(c\in U\), \[\rank[c]=\begin{cases}0&\text{ if }\lambda[\square ,c]=\emptyset\\
  \sup[\rank[\square ]\mid \lambda[\square , c]=\y\,]+1&\text{ if }\lambda[\square ,c]\neq\emptyset\end{cases}.\tag{\(\ast\)}\]
  Further, this function is unique.
\end{lemma}
 Condition (\(\ast\)) implies: elements of rank \(r\) are \(\lambda\)-minimal in the complement of elements of lower rank. To condense the notation a bit, for \(h\in \PP[U]\) define \(\suprank[h]:=\sup[\rank[\square ]\mid \square \in h]\). (\(\ast\)) then becomes 
 \[\rank[c]=\suprank[\lambda[\square ,c]]+1.\]
 \medskip
 
 \noindent Proof of the lemma:\quad First we see that (\(\ast\)) implies the definition. If \(S\subset U\)  then restriction gives a function \(\rank\colon S\to \WW\). If \(S\) is nonempty then the image is nonempty, and since \(\WW\) is well-ordered there is a \(\leq\)-minimal element in the image. Suppose \(\rank[a]\) is minimal. 
If \(\lambda[x,a]=\y\) then the rank condition in (2) implies \(\rank[x]<\rank[a]\), which implies \(x\notin S\). The element \(a\) is therefore minimal.

Now assume (1). To proceed by recursion  we need a formulation in terms of partially-defined functions on \(\WW\). The inverse is a function that takes \(r\in \WW\) to the set of elements of rank \(r\):   \(\rank^{-1}[\square ]\in \PP[U]\). Equation (\(\ast\)) above becomes:
\[f[r]=\begin{cases}\emptyset &\text{ if \(r\) is a limit}\\
\min_{\lambda}[U-(\cup_{\square \leq r-1}f[\square ])] &\text{ if not. }
\end{cases}\]
 We can interpret this as a recursion relation on the collection of pairs \((r,f)\), where \(r\in \WW\) and \(f\colon (\square <r)\to \PP[U]\). Recursion therefore gives a maximal recursive \(F\). The recursion relation implies that the \(F[\square ]\) are disjoint. Maximality implies that every element in \(U\) is in the image. Together these imply we can invert to get a rank function. 
 
 Uniqueness follows from uniqueness of recursion.\qedms
 
The proof shows:
\begin{corollary} An injective morphism of well-founded pairings with transitive image preserves rank functions.
\end{corollary}
The lemma also implies item (3) in the theorem. The point is that while the ZFC axioms include well-founded, ``set'' there means ``set of the form \(\in\!\![\square ,a]\) for some \(a\)'', rather than set in the sense used here. However the full axom scheme enables internal construction of a  rank function in essentially the same way (\cite{Jech}, Example 6.2), which then implies `well-founded' in the stronger sense.
\subsection{Proof of the Theorem}Parts are labeled separately
 \medskip\par\noindent\textbf{Proof of (1)} 
 The plan is to construct an appropriate rank function on \(\WW\), then fill in  to get \(\lambda^{adj}\). The rank will be monotone increasing, so determined by a function \(\minrank\colon \WW\to \WW\). Explicitly, \(b\) has rank \(r\) if \(\minrank[r]\leq b<\minrank[r+1]\). 
 
Now set up for recursion. The domain of the recursion relation is \((a,f)\) with \(a\in \WW\) and \(f\colon(\square <a)\to \WW\) monotone increasing. 
The recursion conditon \(R\colon \dom[R]\to \WW\)\ is:
\begin{enumerate}
 \item \(R[0,f]=0, R[1,f]=1\);
 \item If \(r\) is a limit then \(R[r,f]=\sup[f[\square <r]]\); otherwise
 \item if \(r-1\) is a limit, then \(R[r,f]=f[r-1]\); otherwise
 \item if  \(r-2\) is a limit, then \(R[r,f]=\card\bigl[\PP[\square <r-2]-\bP[\square <r-2]\bigr]\) (see below for `\(\card\)'); and otherwise, finally
 \item \(R[r,f]=f[r-1]+\card\bigl[\PP[\square \leq r-1]-\PP[\square < r-1]\bigr]\).
\end{enumerate}
Recursion then implies there is a maximal recursion function, with domain \(\WW\). Denote this function by \(\minrank[\square ]\)\footnote{This is the iterated-powerset function usually called `Beth'. The recursion condition could be considerably condensed, at the cost of added argument with cardinal arithmetic.}. 
\subsubsection*{Cardinality} Cardinals have a deep and extensive literature. We give the tiny fragment needed for (4) and (5) above, then resume the proof of (1). Suppose \(A\) is a set. We define the \textbf{cardinality}, denoted \(\card[A]\in \WW\) to be the minimum equivalence class  of the form \(\langle A, \leq\rangle\), where `\(\leq\)' is a well-order on \(A\). ``Cardinals'' are elements of \(\WW\) that are cardinalities of some set. It is immediate from this definition that:
\begin{enumerate}
 \item \(\card[A]=\card[B]\) if and only if there is a bijection \(A\simeq B\);
 \item if \(c\) is a cardinal then \(\card[\square <c]=c\); and
more generally  \item  if \(r\in \WW\) and \(c\) is a cardinal then \(\card[r\leq \square \leq r+c]=c\).
\end{enumerate}
The reason cardinality is used in the recurrence relation is that such relations must be well-defined, and cardinality gives a well-defined measure of the size of a set. 

We return to the proof of item (1) of the Theorem. Use the \(\minrank\) function found above  to define a rank function on  \(\WW\). The recursion relation unwinds to show: for each \(r\in \WW\) there exists a bijection from elements of rank \(r\) to functions in \(\bP[\WW]\) of suprank \(r-1\). Note that (3) implies that if \(r\) is a limit, there are no elements of rank \(r\), so \(r-1\) is defined. Using the axiom of Choice we can choose such an isomorphism \(E_r\) for each \(r\).  Putting these together we get a bijection \(E\colon\WW\to\bP[\WW]\). Define the pairing \(\epsilon\colon \WW\times\WW\to \yn\) to be the adjoint of this. By construction \(\epsilon\) satisfies condition \((\ast)\) of 4.1, so is well-founded. This concludes the proof of (1).\qedms

\medskip\par\noindent\textbf{Proof of (2)}.\quad Suppose \(\epsilon\) is a well-founded pairing on \(\WW\) with bijective adjoint, and \((U,\lambda)\) is some other well-founded pairing. The claim is that there is a unique injective morphism of pairings \(\lambda\to \epsilon\) with transitive image. Denote rank functions for these by \(\rank_{\epsilon}, \rank_{\lambda}\). Since the desired morphism would have to preserve rank, we show its existence by recursion on rank. 

Define the domain of the recursion relation to be \((r,f)\) where \(r\in \WW\) and \(f\colon(\square \in U\mid \rank_{\lambda}[\square ]<r)\to 
(\square \in \WW\mid \rank_{\epsilon}[\square ]<r)\) a rank-preserving injection with transitive image. Define \(R[r,f]\colon ( \rank_{\lambda}[\square ]=r)\to (\rank_{\epsilon}[\square ]=r)\) by the composition indicated in the diagram
  
\[\xymatrix{  ( \rank_{\lambda}[\square ]=r)\ar[d]^{\lambda^{adj}}\ar@{ -->}[r]^{R[r,f]}& ( \rank_{\epsilon}[\square ]=r)\\
 \PP[ \rank_{\lambda}[\square ]<r]\ar[r]^{f_{!}}&\PP[ \rank_{\epsilon}[\square ]=r]\ar[u]^{\simeq}_{(\epsilon^{adj})^{-1}}
 }.\]
 Here \(f_{!}\) is the push-forward defined, for \(h\in \PP[ \rank_{\lambda}[\square ]<r]\), by
 \[(f_{!}[h])[\square ]=
\begin{cases}h[f^{-1}[\square ]]&\text{ for }\square \in\im[f]\\
 \n &\text{ for }\square \notin\im[f].
\end{cases}\]
This uses: \(f\) is defined and injective on \((\rank_{\lambda}[\square ]<r])\); and the image is the support of a 2-valued function, denoted here by `\((\square \in\im[f])\)'.

We claim the result of recursion with this relation has the desired properties. Injectivity follows from injectivity of \(\lambda^{adj}\); it is a morphism of pairings because \(f^*[f_{!}[h]]=h\); and \(f_{!}[h]\) vanishing off the image of \(f\) implies that the image is transitive. \qedms
\medskip\par\noindent\textbf{Proof of (4)}.\quad
The claim is that the universal pairing of (1) satisfies the ZFC-1 axioms. We expand on the ``-1'': In ZFC, `functions' are defined via sets: a primitive function \(f\colon A\to B\) is a ZFC function if its graph is a ZFC subset of \(A\times B\). Having enough functions for the basic business of set theory therefore depends on finding sufficiently many such subsets. It was a discovery of Zermillo that admitting first-order logical formulas in the other ingredients does this. In the approach here, `functions' are primitives and first-order logic is not needed to identify which ones should be included in the theory.  The `-1' in `ZFC-1' means that the relevant axioms are satisfied for \emph{all} primitive functions, not just those given by first-order formulas. 

The ZFC axioms, as listed in \cite{Jech} Ch. 1, and the reasons \(\epsilon\) satisfies them, are:
 \begin{description}
\item[Extentionality] Sets are determined by their elements. This is injectivity of the adjoint;
\item[Union]  \(\exists (\cup a)\in  U \) such that \(\epsilon[\square ,\cup a]=(\exists b\mid \epsilon[\square ,b]\ad \epsilon[b,a])\). In other words, the expression on the right is in the image of \(\epsilon^{adj}\). This is so because  \(\epsilon^{adj}\) is a bijection onto  \(\bP[\WW]\).
\item[Powerset]  \(\exists (\PP[a])\in  U \) such that for all \(b\), \(\epsilon[b,\PP[a]]=(\epsilon[\square ,b]\subset \epsilon[\square ,a])\). Again this is surjectivity of \(\epsilon^{adj}\) onto  \(\bP[\WW]\).;
\item[Infinity] There is \(a\) with \(\epsilon[\square ,a]\) infinite. This is equivalent to our primitive hypothesis that the natural numbers is a set;
\item[Choice] This is a primitive hypothesis here;
\item[Separation] If \(P\) is a 2-valued function on \( U \) given by a first-order formula in the set operations then the intersection of a set with \(P\) is also a set. A set is bounded in \(\WW\), so its intersection with \(P\) is also bounded. But here `bounded' implies `set',  whether or not the function is given by a formula;
\item[Well-founded, or Regular] \(\epsilon\) is  well-founded by construction; 
\item[Replacement] Suppose \(A\) is a set and \(f\colon  A \to  U \) is a ZFC function, then the image of \(f\)  is  a set. Here the image of a set under any function, ZFC or not, is a primitive set by \ref{sssect:subobjects2}. 
\end{description}
This completes the proof of (4), and therefore the theorem.\qed

\section{Coherent limits} \label{sect:coherent-lims}
 In the previous section we saw that the set theory defined here corresponds to a maximal model for ZFC. Here we find a property that characterizes this maximal model. Models that do not satisfy the condition are inconsistent with mainstream practice, and some of them are not even consistent with standard presentations of elementary calculus.

\subsection{Models and universes}
A  \textbf{model} is an implementation of the logic that satisfies the axioms. This consists of: a Universe  \( U \); an equality pairing (so it is a pre-set);  and a 2-valued ``membership'' pairing  \(\in\colon  U \times  U \to  \yn\). 
 
Traditional accounts  of set theory (eg.~\cite{Jech})  omit references to a universe. The basics of set theory date back to the late 1800s, and this is one of the ways it remains pre-modern.  There are also two points of view: One is that \(U\) is a universe of sets, and  \(a\in b\) makes sense because ``an element of a set is again a set''.  A mathematically more-coherent version is that \(U\) is a universe of possible elements, and sets are collections of elements specified by 2-valued functions of the form \((\square \in a)\), for some \(a\in U\). In other words `\(\in\)' gives a \emph{parameterization} of sets by elements, not an identification. 
\subsubsection{Sets are sets}\label{sssect:setsRsets}
The logic of axiomatic set theory has equality and quantification built in, so sets in such a theory should automatically be sets in the sense used here. There is, however, a loophole. Quantification here requires being able to identify the empty function among all 2-valued functions. Quantification in a model of  ZFC only requires being able to identify it among 2-valued functions \emph{in the model}, and so as far as we know now, \(\WW\) might be a set in a very thin model. This seems unlikely, but if it does happen then \(\WW\) is a regular strong-limit cardinal in the theory, and we truncate it there to get a ZFC model in which the sets are sets.

 \subsection{Coherent Limit Property} \label{ssect:coherent-limit} The formulation is chosen because it concerns a key step in recursion.
 
 A model  \((U,\in)\)  for the ZFC\(\pm1\) axioms has the \textbf{coherent limit property} if the following is always true:
 
\textbf{If} \(A \) is a well-ordered \(U\)-set,  and 
 for each \(a\in A\) a 2-valued \(U\)-function \(f_a\colon (\square < a)\to \yn\)  is given, such that 
 if \(b\leq a\) then the restriction \(f_a\upharpoonright (\square < b)\) is \(f_b\), \textbf{then}
 the formula \[F[a]:=f_c[a]\text{ for any }c>a\] defines a \(U\) function \(A\to \yn\).\medskip

 \begin{proposition} \((U,\in)\) has the Coherent Limit Property  if and only if it is isomorphic to (a truncation of) the maximal model \((\WW,\epsilon)\).\end{proposition}
Here a ``truncation'' is the collection of sets with cardinality less than \(k\), where \(k\) is a regular strong-limit cardinal. As is  well known, this is also a ZFC set theory. 

A recursion relation (\S \ref{ssect:recursion}) generates such a family \(f_a\). The output of recursion is the function given by the formula. But   if the model is not maximal, the Proposition shows that additional argument is needed to ensure the function is in the model. 
The treatment of recursion in  \cite{Jech} p.~22  seems to gloss over this, so as it stands it seems to be  incomplete. Such additional arguments are not part of mainstream methodology, so this methodology implicitly takes place in the maximal theory. For the same reason, standard presentations of elementary calculus are only valid in models that have the same real numbers as the maximal model. 

\subsection{Coherent unions}
If we specialize the statement to 2-valued functions (\(B=\yn\)) then we get a reformulation in terms of subsets: a model has the \textbf{coherent union property} if the following is true:

\textbf{If} \(A \) is a well-ordered \(U\)-set,  and 
 for each \(a\in A\) a \(U\) subset \(S_a\subset (\square <a)\)  is given, such that 
 if \(b\leq a\) then the intersection \(S_a\cap (\square < b)\) is \(S_b\), \textbf{then}
 the union \(\bigcup_{\square \in A}S_{\square }\) is  a \(U\) subset of \(A\).\medskip

This is equivalent to the coherent limit property, and again is only true if \((U,\in)\) is maximal. In this form it is easier to explain the problem. The ZFC Union axiom asserts that the union is a \(U\) subset of \(A\) \emph{if} the collection \(\{S_{\square }\}\subset \PP_U[A]\) is a \(U\) subset. But the hypotheses of the property do not imply \(\{S_{\square }\}\) is a \(U\) subset, and if \(U\) is not maximal  then the Proposition asserts that we can find examples where the conclusion of the property fails. 

\subsection{Proof of the Proposition} We prove the union form.
If \((U,\in)\) is a truncation of \((\WW,\epsilon)\) then the property is a very special case of the union lemma. 

For the converse, ``not a truncation'' means there is a \(U\)-set \(A\) that contains a non-\(U\) set in the sense defined here.  Equivalently, there is an \(\WW\) function \(F\colon A\to \{0,1\}\) that is not a  \(U\) function. 
We may assume that  \(A\) is well-ordered, and is minimal in the sense that if \(a\in A\) then every primitive function  \((\square < a)\to \{0,1\}\) is a \(U\) function. For each \(a\in A\) define \(f_a\) to be the restriction of \(F\) to \((\square <a)\). Minimality implies these functions satisfy the hypotheses of the Property, but the limit does not satisfy the conclusion. \qed

\section{Appendix:  Sets without the Quantification hypothesis}\label{sect:noQuantification} Of the primitive hypotheses in \S\ref{ssect:primitives}, Two and Infinity are clearly necessary, and Choice has been extremely well-tested. 
 The justification for the Quantification Hypothesis (QH) is less clear, so in this appendix we describe the changes needed to avoid it. It turns out the set theory is essentially the same as with QH, and gives no hint QH is in any way problematic. 
\subsection{Outcomes}
First, some notation. We say that a pre-set is \(Q^1\) if it supports quantification, \S4.3. A pre-set \(A\) is \(Q^n\) for \(n>1\) if it is \(Q^1\) and \(\PP[A]\) is \(Q^{n-1}\). \textbf{Sets} are \(Q^{\infty}\) pre-sets. The Quantification hypothesis is that \(Q^1\Rightarrow Q^{\infty}\), so in this case the intermediate notations are unnessary. If QH fails then there are \(Q^1\) pre-sets that are not \(Q^{\infty}\).
Further, there is a maximal strong-limit cardinal we denote by \(q_{\infty}\). \(Q^{\infty}\) pre-sets are those with cardinality less than this, and there are some  finite-quantification pre-sets with larger cardinality. This is a version of von Neumann's ``Axiom of size''. With QH there is no useful analog of this Axiom: a pre-set must be \(Q^1\) for cardinality to be defined, so is already a set. 

As in \S\ref{ssect:w-oClassification},  \(\WW\) denotes the order-equivalence classes of well-ordered \(Q^1\) pre-sets.  Let \(Q^j\subset \WW\) denote the transitive sub-presets that are \(Q^j\), for \(j\geq 1\). Comparing with the definition above we see \(Q^{\infty}\)=\((\square <q_{\infty})\).

The  construction in Section 3 still applies to give a pairing on \(Q^{\infty}\) that is universal for well-founded pairings on sets. As in the QH case this pairing  satisfies the ZFC-1 axioms  \emph{except} for the Union axiom. The Union lemma  only gives \(Q^1\). QH promotes this to \(Q^{\infty}\), but not[QH] does not. 

The ZFC Union axiom is equivalent to regularity of the maximal strong-limit cardinal.  Success with the Union axiom in traditional set theory  supports adding this as a hypothesis, but we see no reason not to assume QH and get Union for free. 

 In not[QH] the same argument as in Section \ref{sect:wfpairings} gives a universal well-founded pairing \(\epsilon\colon Q^{\infty}\times Q^{\infty} \to \yn\) with all the same properties except, as noted above, the   Union axiom. 
The  observations about the Coherent Limit Property apply whether or not QH is assumed.


\enddocument